\documentclass[11pt ]{article}
\usepackage[pdftex]{graphicx}

 \newtheorem{Thm}{Theorem.}
\newtheorem{Prop}{Proposition.}
\newtheorem{Cor}{Corollary.}
\newtheorem{Lem}{Lemma.}

\begin{document}

\centerline {\bf A new look at the Julg-Valette homotopy for groups acting on trees.}

\bigskip
\centerline { \bf Pierre Julg}

\centerline { MAPMO (Universit\'e d'Orl\'eans et CNRS)}

\bigskip\bigskip\bigskip

In [JV1] [JV2] we proved the K-theoretic amenability of localy compact groups acting on trees (with amenable stabilizers) by combining two ingredients: 

1)  We construct a very simple Fredhom module associated to the group action on a tree. It defines a class $\gamma$ in Kasparov's ring $KK_G({\bf C},{\bf C})$.

2) We show that $\gamma=1$ by constructing a homotopy using the simple but  non trivial fact that the distance kernel on the set of vertices of the tree is of conditionnallly negative type. Or equivalently the existence of an affine action of $G$ on the $\ell^2$ space of edges of the tree.

Note that the two ingredients are of different nature. The first can be generalized to other situations such as Bruhat-Tits buildings [JV3][KS1], or (hyper)bolic spaces [KS2]. The second is very specific to trees, or generalizations (e.g. CAT(0) cubic complexes cf [B]).

We present here a perhaps more natural proof of the same result. It is inspired by the construction by Pytlik and Szwarc [PS] of a family of uniformly bounded representations of a free group, generalized by Valette [V1] and Szwarc [S1] to groups acting on trees. We have had this new version for quite a long time in our private notes. We thank Amaury Freslon and Jacek Brodzki for convincing us that making these notes available could be useful to others. We thank Professor Szwarc for pointing us reference [S1]. We dedicate this short paper to the memory of Tadeusz Pytlik, who died in 2006 [S2].

\bigskip\bigskip

{\bf 1. Notations:}
\bigskip

Let $X=(X^0, X^1)$ be a tree.  There is no orientation on $X$ so that the set $X^1$ of edges is just a subset of $X^0\times X^0$ stable by the map $(x,y)\mapsto (y,x)$. By hypothesis, for any $x$ and $y$ in $X^0$ there is a unique path joining $x$ to $y$. We denote by $[x,y]$ the set of vertices lying between $x$ and $y$ and $d(x,y)$ the number of edges between $x$ and $y$.
\bigskip
Let $G$ be a group acting on $X$. In other words $G$ acts on $X^0$ in such a way that the subset $X^1\subset X^0\times X^0$ is stable.
\bigskip
We consider the Hilbert space $\ell^2(X^0)$, $(\delta_x \ ,x\in X^0)$ its canonical Hilbert basis and $\pi_0$ the unitary representation arising from the action of $G$ on the set $X^0$ defined by $\pi_0(g)\delta_x=\delta_{gx}$. 
Let $\ell^2(X^1)^-$ be the quotient  of $\ell^2(X^1)$ by the subspace generated by the vectors $\delta_{(x,y)}+\delta_{(y,x)}$ for $(x,y)\in X^1$ and $\pi_1$ the unitary representation of $G$ on $\ell^2(X^1)^-$ defined by 
$\pi_1(g)\delta_{(x,y)}=\delta_{(gx,gy)}$

\bigskip

For $x\in X^0$ let $V(x)$ be the set of neighbours of $x$. Let the cardinality of $V(x)$ be denoted $q_x+1$. We assume for simplicity that $q_x$ is bounded.
We define the following bounded operators $S$ and $Q$ on $\ell^2(X^0)$ by the following formulae:

 $$S\delta_x=\sum_{y\in V(x)}\delta_y$$
 $$Q\delta_x=q_x\delta_x$$
 
 Note that the operators $Q$ and $S$ only depend on the tree structure on $X$ and therefore commute with the unitaries $\pi_0(g)$.
\bigskip\bigskip

{\bf 2. The Pytlik-Szwarc operator.}

\bigskip
We now choose an origin $x_0\in X^0$. Let $p_0$ the orthogonal projection onto the vector $\delta_{x_0}$.
Let us define the operator $P$ on $\ell^2(X^0)$  by
$$P\delta_{x_0}=0$$
$$P\delta_x=\delta_{x'}$$
where for $x\neq x_0$ we define $x'$ to be the unique neighbour of $x$ lying between $x_0$ and $x$.
\bigskip

\begin{Prop}
The operator $P$ such defined is bounded on $\ell^2(X^0)$ and satisfies:
$$PP^*=Q+p_0$$
$$P+P^*=S$$
\end{Prop}

Proof: For $ x\neq x_0$, one has $P^*\delta_x=\sum_{y\in V(x)\setminus \{x'\}}\delta_y$ so that clearly $(P+P^*)\delta_x=S\delta_x$. Similarly  $(P+P^*)\delta_{x_0}=P^*\delta_{x_0}=\sum_{y\in V(x_0)}\delta_y=S\delta_{x_0}$.

On the other hand, for $x\neq x_0$ and $y\in V(x)\setminus \{x'\}$ one has $P\delta_y=\delta_x$ so that $PP^*\delta_x=q_x\delta_x$,  
whereas for any  $y\in V(x_0)$,  $P\delta_y=\delta_{x_0}$ so that
 $PP^*\delta_{x_0}=(q_{x_0}+1)\delta_{x_0}$

\bigskip
\begin{Cor}
Let $T_t=1-tP+((1-t^2)^{1/2}-1)p_0$ for any  $t\in [0,1]$
Then the operator $T_t$ satisfies the following formula:
$T_tT^*_t=1-tS+t^2Q$. In particular, the operator $T_tT^*_t$
and commutes with the unitaries $\pi_0(g)$, $g\in G$.
\end{Cor}

Indeed we have $T_t=1-tP+\alpha p_0$ where $\alpha$ satisfies $\alpha^2+2\alpha +t^2=0$.
A straightforward calculation (using the obvious fact that $Pp_0=0$) yields $T_tT_t^*=1-t(P+P^*)+t^2PP^*+(2\alpha+\alpha^2)p_0$ and the result follows from the proposition.
\bigskip\bigskip

{\bf 3. Construction of new representations.}

\bigskip
Let us consider the space $D(X^0)$ of finitely supported functions on $X^0$ as a dense subspace of $\ell^2(X^0)$. Clearly $D(X^0)$ is stable by the operators  $P$, $P^*$ and $\pi_0(g)$ for $g\in G$.

\bigskip
\begin{Lem}
 For any complex number $z$ the operator $$(1-zP)^{-1}=\sum_{k=0}^{\infty}z^kP^k$$ is defined on $D(X^0)$ and one has $$(1-zP)^{-1}\delta_x=\sum_{y\in [x_0,x]}  z^{d(y,x)}\delta_y$$
\end{Lem}

Proof. Let the elements of $[x_0, x] $ be denoted $x_0$, $x_1$,..., $x_n=x$. Then $P^k\delta_x=\delta_{x_{n-k}}$ for $k\leq n$ and $0$ for $k>n$. Note that $d(x_{n-k},x)=k$. This makes the statement straightforward.
\bigskip

\begin{Thm}
 (Pytlik-Szwarc)  Let $z$ be a complex number such that $\vert z\vert <1$. For any $g\in G$ the operator $\rho_z(g)=(1-zP)^{-1}\pi_0(g)(1-zP)$ extends to a bounded operator on $\ell^2(X^0)$, thus defining a representation $\rho_z$ of $G$ in $\ell^2(X^0)$. The operator  $\rho_z(g)-\pi_0(g)$ is a finite rank operator and the representation $\rho_z$ is uniformly bounded, i.e. 
$$\sup_{g\in G}\Vert \rho_z(g)\Vert <\infty$$
\end{Thm}
\bigskip

\begin{Thm}
  For any real number such that $0<t<1$  the operator $\tilde\rho_t(g)=T_t^{-1}\pi_0(g)T_t$ extends to a {\it unitary} operator on $\ell^2(X^0)$, thus defining a unitary representation $\tilde\rho_t$ of $G$ in $\ell^2(X^0)$. The operator $\tilde\rho_t(g)-\pi_0(g)$ is a finite rank operator. The uniformly bounded  representation $\rho_t$ is equivalent to the unitary representation $\tilde\rho_t$.
\end{Thm}

Proof: Let us prove that $\rho_z(g)-\pi_0(g)$ has finite rank and a norm bounded independently from $g$. It is enough to consider the operator $$\rho_z(g)\pi_0(g)^{-1}-1=z(1-zP)^{-1}(P-P')$$ where $P'$ is defined just as $P$ but replacing $x_0$ by $gx_0$. It is clear that the above operator is supported on the finite dimensional subspace generated by the $\delta_x$ for $x\in [x_0,gx_0]$. That subspace is indeed stable by $P$ and $P'$ which have restrictions of norm 1. Therefore the norm of $(1-zP)^{-1}(P-P')$ is at most $2\sum \vert z\vert^k=2(1-\vert z\vert )^{-1}$. This proves theorem 1.

To deduce theorem 2 note that $\tilde \rho_t(g)=u_t^{-1}\rho_t(g)u_t$ where $u_t=(1-p_0)+(1-t^2)^{1/2}p_0$ is an invertible operator which differs from the identity by a compact operator. It remains to show that $\tilde\rho_z(g)$ is unitary. Since it is invertible it is enough to compute $\rho_t(g)^*\rho_t(g)=T_t^*\pi_0(g)^{-1} {(T_tT_t^*)}^{-1}\pi_0(g)T_t$ which is equal to 1 since $T_tT_t^*$ commutes to $\pi_0(g)$ by the corollary to the proposition above.
\bigskip

{\bf Remark. }
The link with the original Julg-Valette approach is given by an easy computation: the kernel of the (densely defined) operator $(T_tT_t^*)^{-1}$ is 
$$<T_t^{-1}\delta_x ,T_t^{-1}\delta_y>= t^{d(x,y)}=e^{-\lambda d(x,y)}$$ if $t=e^{-\lambda}$.

\bigskip\bigskip
{\bf 4. The limit when $t$ tends to 1.}
\bigskip

Let us now calculate the limit of $ \tilde \rho_t(g)$ when $t\rightarrow 1$. Let us recall the definition of the Julg-Valette map $F: \ell^2(X^0)\rightarrow \ell^2(X^1)^-$: 
 $$F\delta_{x_0}=0$$
 $$F\delta_x=\delta_{(x,x')}$$
 
 We have $Fp_0=0$, $F^*F=1-p_0$ and $FF^*=1$.
 \bigskip
 
 \begin{Lem}
  Let $b:\ell^2(X^1)^-\rightarrow \ell^2(X^0)$ be the coboundary operator defined by $b\delta_{(x,y)}=\delta_y-\delta_x$. Then the operators $F$ and $P$ are related by the formulae:
 $$1-P=bF+p_0$$
 $$(1-P)F^*=b$$
 \end{Lem}

Indeed, $(1-P)\delta_x=\delta_x-\delta_{x'}=bF\delta_x$ if $x\neq x_0$ and $(1-P)\delta_{x_0}=\delta_{x_0}$. The  second formula follows from the first since $(1-P)F^*=bFF^*=b$. 
\bigskip

{\bf Remark.} Let $c(x,y)=\sum\delta_{(x_i,x_{i+1})} $
if the elements of $[x, y] $ are denoted $x_0=x$, $x_1$,..., $x_n=x$. This is the cocycle realizing explicitely the Haagerup property for groups acting properly on trees: $\Vert c(x,y)\Vert^2=d(x,y)$. 
On has $bc(x,y)=\delta_y-\delta_x$, hence by the second formula above,   $c(x,y)=F(1-P)^{-1}(\delta_y-\delta_x)$.

\bigskip
\begin{Cor}
 The operator $(1-P)^{-1}b$ extends to a bounded operator and one has:
$$(1-P)^{-1}b=F^*$$
\end{Cor}

It follows indeed from the lemma that $(1-P)F^*=b$.

\bigskip
\begin{Prop}
 For any $g\in G$ the unitary operator $\tilde \rho_t(g)$ ($0<t<1$) converges strongly to
$$\tilde \rho_1 (g)=F^*\pi_1(g)F+p_0$$
when $t\rightarrow 1$.
\end{Prop}

Proof: Let us first prove that $F\tilde \rho_t(g)F^*$ strongly converges to $\pi_1(g)$. 
One has 
$$F\tilde \rho_t(g)F^*=FT_t^{-1}\pi_0(g)T_tF^*=F(1-tP)^{-1}\pi_0(g)(1-tP)F^*$$
 which evaluated on $\delta_e$ ($e\in X^1$) converges to 
$F(1-P)^{-1}\pi_0(g)(1-P)F^*$  i.e. 
$$F(1-P)^{-1}\pi_0(g)bFF^*=F(1-P)^{-1}b\pi_1(g)=FF^*\pi_1(g)=\pi_1(g).$$

We deduce that $(1-p_0)\tilde \rho_t(g)(1-p_0)$ converges strongly to $F^*\pi_1(g)F$.

On the other hand we check that  $\tilde \rho_t(g)\delta_{x_0}\rightarrow\delta_{x_0}$ when $t\rightarrow 1$. We have indeed $T_t\delta_{x_0}=(1-t^2)^{1/2}\delta_{x_0}$ so that 
$(1-tP)^{-1}\pi_0(g)T_t\delta_{x_0}=(1-t^2)^{1/2}(1-tP)^{-1}\delta_{gx_0}=(1-t^2)^{1/2}\sum t^{d(y,gx_0)}\delta_y$ where the sum is extended to the vertices $y$ of the segment $[x_0,gx_0]$.
Finally $T_t^{-1}\pi_0(g)T_t\delta_{x_0}=t^{d(x_0,gx_0)}+(1-t^2)^{1/2}\sum t^{d(y,gx_0)}\delta_y$ where $x_0$ is now excluded from the sum. Hence the result.

As a consequence $\tilde \rho_t(g)p_0$  converges normally to $p_0$, and since $\tilde \rho_t$ is unitary, replacing $g$ by $g^{-1}$ we also have that $p_0\tilde \rho_t(g)$ converges normally to $p_0$. 

The proposition clearly follows.

\bigskip\bigskip
{\bf 5. Classes in $KK_G$-theory.}

\bigskip
Let us consider the Hilbert space $\ell^2(X^0)$ equipped with the two representations $\tilde \rho_t$ and $\pi_0$,  which differ by  compact operators. The triple $(\tilde \rho_t ,\pi_0, Id)$ defines an element of the Kasparov group $KK_G({\bf C},{\bf C})$. It is independent of the value on $t\in[0,1]$ since the operators $\tilde \rho_t(g)$ are  strongly continuous in $t$. Now when $t=0$ we have $\tilde \rho_0=\pi_0$ so that the element is equal to $0$. On the other hand when $t=1$ we have $\tilde \rho_1 (g)=F^*\pi_1(g)F+p_0$ so that the element  is equal to $1-\gamma$ where $\gamma$ is defined as in Julg-Valette by the triple $(\pi_0, \pi_1,F)$.
\bigskip
\begin{Cor}
 We have the equality $\gamma=1$ in the Kasparov group $KK_G({\bf C},{\bf C})$.
\end{Cor}

\bigskip\bigskip

{\bf Bibliography:}

\bigskip
[B] Jacek Brodzki, {\it A differential complex for groups acting on CAT(0) cube complex}.

[JV1] Pierre Julg, Alain Valette,  {\it K-moyennabilit\'e pour les groupes op\'erant sur les arbres}, C. R. Acad. Sci. Paris S\'er. I Math. 296 (1983), no. 23, 977Ð980. 

[JV2] Pierre Julg, Alain Valette,{\it  $K$-theoretic amenability for $SL_2({\bf Q}_p)$, and the action on the associated tree},  J. Funct. Anal. 58 (1984), no. 2, 194Ð215.

[JV3] Pierre Julg, Alain Valette, {\it Fredholm modules associated to Bruhat-Tits buildings}. Miniconferences on harmonic analysis and operator algebras (Canberra, 1987), 143Ð155, Proc. Centre Math. Anal. Austral. Nat. Univ., 16, Austral. Nat. Univ., Canberra, 1988.

[KS1]Gennadi Kasparov, Georges Skandalis, {\it  Groups acting on buildings, operator K-theory, and Novikov's conjecture}, K-Theory 4 (1991), no. 4, 303Ð337.

[KS2] Gennadi Kasparov, Georges Skandalis, {\it  Groups acting properly on "bolic'' spaces and the Novikov conjecture}. Ann. of Math. (2) 158 (2003), no. 1, 165Ð206. 
 
 [PS] Tadeusz Pytlik,  Ryszard Szwarc,  {\it An analytic family of uniformly bounded representations of free groups},  Acta Math. 157 (1986), no. 3-4, 287Ð309.
 
 [S1] Ryszard Szwarc, {\it Groups acting on trees and approximation properties of the Fourier algebra}, J. Funct. Anal 95 (1991), 320-343.
 
 [S2] Ryszard Szwarc, {\it Tadek Pytlik in my memories}, Colloquium Mathematicum, 108 (2007).
 
 [V1] Alain Valette, {\it Cocycles d'arbres et repr\'esentations uniform\'ement born\'ees.}  C. R. Acad. Sci. Paris S\'er. I Math. 310 (1990), no. 10, 703Ð708. 
 
  [V2] Alain Valette, {\it  Les repr\'esentations uniform\'ement born\'ees associ\'ees ˆ un arbre r\'eel},  Bull. Soc. Math. Belg. S\'er. A 42 (1990), no. 3, 747Ð760. 

[V3] Alain  Valette,  {\it Addendum: "Uniformly bounded representations associated with a real tree''.}  Bull. Soc. Math. Belg. S\'er. A 44 (1992), no. 1, 101Ð102.

  \end{document}